\documentclass{article}
\usepackage[a4paper, total={7in, 10in}]{geometry}
\usepackage[T1]{fontenc}
\usepackage{graphicx}
\usepackage{mathtools}
\usepackage{amsmath,amsthm,amssymb,mathrsfs,amstext, titlesec,enumitem, comment, graphicx, color, xcolor, stmaryrd,mathabx}
\usepackage{thmtools}
\usepackage{xpatch}
\usepackage{tikz-cd}
\usepackage{nameref}
\usepackage{hyperref}
\makeatletter
\def\Ddots{\mathinner{\mkern1mu\raise\p@
\vbox{\kern7\p@\hbox{.}}\mkern2mu
\raise4\p@\hbox{.}\mkern2mu\raise7\p@\hbox{.}\mkern1mu}}
\makeatother
\newtheorem{theorem}{Theorem}[section]

\newtheorem{lemma}[theorem]{Lemma}

\theoremstyle{definition}

\newtheorem{ex}{Example}
\newtheorem{definition}[theorem]{Definition}

\begin{document}
\title{\textbf{The Interplay between Additive and Multiplicative Central Sets Theorems } }

	\date{}
	\author{Pintu Debnath
		\footnote{Department of Mathematics,
			Basirhat College,
			Basirhat-743412, North 24th parganas, West Bengal, India.\hfill\break
			{\tt pintumath1989@gmail.com}}
		\and
		 Sayan Goswami 
		\footnote{Department  Mathematics,
			Ramakrishna Mission Vivekananda Educational and Research Institute,
			Belur, Howrah, 711202, India\hfill\break
			{\tt sayan92m@gmail.com}	
	}
    \and 
    Chunlin Liu
    \footnote{School of Mathematical Sciences, Dalian University of Technology, Dalian, 116024, P.R. China  
   }
      \footnote{Institute of Mathematics, Polish Academy of Sciences, ul. Śniadeckich 8, 00-656 Warszawa, Poland \hfill\break 
    	{\tt chunlinliu@mail.ustc.edu.cn}}
}
\maketitle	

 \begin{abstract}
The concept of Central sets, introduced by Furstenberg through the framework of topological dynamics, has played a pivotal role in combinatorial number theory. Furstenberg's Central Sets Theorem highlighted their rich combinatorial structure. Later in \cite{DHS08}, De, Hindman, and Strauss strengthen this theorem using the algebraic framework of the Stone--\v{C}ech compactification. In this article, we establish a unified version of the Central Sets Theorem that simultaneously captures both additive and multiplicative structures.
\end{abstract}


\section{Introduction}

Let \( X \) be a compact Hausdorff space and \( T: X \to X \) a continuous self-map. The pair \( (X, T) \) is what we call a \textbf{dynamical system}. Given an open set \( U \subseteq X \) and a point \( x \in X \), the set  
\[
R(x, U) = \{ n \in \mathbb{N} \mid T^n x \in U \}
\]  
is known as the \textbf{return times set} of \( x \) to \( U \).  

Although this definition is presented in the context of a \(\mathbb{Z}\)-action, it naturally extends to more general dynamical systems that arise from actions of arbitrary semigroups.

In his work \cite[Proposition 8.21]{F81}, Furstenberg introduced the concept of \textbf{central sets} via return times and proved the celebrated \textbf{Central Sets Theorem}, which has since become a cornerstone in topological dynamics and combinatorics.

Before proceeding further, we adopt the following conventions:  
\begin{itemize}
    \item The set \( \mathbb{N} \) denotes the set of all positive integers.  
    \item For any nonempty set \( X \), let \( \mathcal{P}_f(X) \) denote the collection of all nonempty finite subsets of \( X \).  
    \item The notation \( {}^\mathbb{N} X \) represents the set of all sequences over \( X \).  
\end{itemize}

Furstenberg's concept of central sets plays a crucial role in combinatorial number theory and topological dynamics. These sets exhibit rich combinatorial structures and possess strong recurrence properties. One of the fundamental results in this area is the \textbf{Central Sets Theorem}, which provides a deep insight into the additive structure of central sets. We now state this important theorem.

\begin{theorem}   [\textbf{Central Sets Theorem}]\cite[Theorem 8.21]{F81}
\label{cst}
Let \( l \in \mathbb{N} \), and let \( A \subseteq \mathbb{N} \) be a central set. For each \( i \in \{1,2,\dots,l\} \), let \( \langle x_{i,m} \rangle_{m=1}^\infty \) be a sequence in \( \mathbb{N} \). Then there exist a sequence \( \langle b_m \rangle_{m=1}^\infty \) in \( \mathbb{N} \) and a sequence \( \langle K_m \rangle_{m=1}^\infty \) in \( \mathcal{P}_f(\mathbb{N}) \) such that:  
\begin{enumerate}  
\item For all \( m \), \( \max K_m < \min K_{m+1} \).  
\item For every \( i \in \{1,2,\dots,l\} \) and any \( H \in \mathcal{P}_f(\mathbb{N}) \),

\[
\sum_{m \in H} \left( b_m + \sum_{t \in K_m} x_{i,t} \right) \in A.
\]

\end{enumerate}  
\end{theorem}

One can see the articles \cite{BJM17, HM99} for the polynomial extension of the above theorem and its several combinatorial applications.
The \textbf{Stronger Central Sets Theorem} extends the classical \textbf{Central Sets Theorem}, offering a deeper structural understanding of central sets. This enhanced version reveals richer combinatorial patterns and a more refined algebraic framework, making it especially valuable in areas such as Ramsey theory and ultrafilter-based combinatorics. It not only guarantees the existence of certain combinatorial configurations within central sets but also provides a more systematic method for their construction.

Using the Stone--\v{C}ech compactification of discrete semigroups, De, Hindman, and Strauss presented a stronger version of the Central Sets Theorem in \cite{DHS08}, which we now state.

\begin{theorem}[\textbf{Stronger Central Sets Theorem}]\cite[Theorem 2.2]{DHS08}
\label{gcst} 
Let \( (S, +) \) be a commutative semigroup, \( \tau = {}^{\mathbb{N}}S \), and let \( C \subseteq S \) be a central set. Then there exist functions  
\[
\alpha: \mathcal{P}_f(\tau) \to S \quad \text{and} \quad H: \mathcal{P}_f(\tau) \to \mathcal{P}_f(\mathbb{N})
\]  
such that:  

\begin{enumerate}
    \item \label{1.41} For all \( F, G \in \mathcal{P}_f(\tau) \) with \( F \subsetneq G \),  
    \[
    \max H(F) < \min H(G).
    \]
    
    \item \label{1.42} For any \( m \in \mathbb{N} \), and any sequence  
    \[
    G_1 \subsetneq G_2 \subsetneq \cdots \subsetneq G_m \quad \text{in } \mathcal{P}_f(\tau),
    \]  
    with \( f_i \in G_i \) for each \( i = 1, 2, \dots, m \), we have  
    \[
    \sum_{i=1}^m \left( \alpha(G_i) + \sum_{t \in H(G_i)} f_i(t) \right) \in C.
    \]
\end{enumerate}  
\end{theorem}

For sets satisfying the above Central Sets Theorem, one can see \cite{J15, L12}.
For a more general version of the above theorem, one can see the article \cite[Theorem 3.4]{DHS09}.  For a simple noncommutative version, we refer to the article \cite{J}.
In this article, our primary goal is to prove the following version of the Central Sets Theorem, which unifies both the additive and multiplicative Central Sets Theorems, that means it unifies the Central Sets Theorem on both $(\mathbb{N},+)$ and $(\mathbb{N},\cdot ).$ The core idea behind this result arises from the following observation: for any finite coloring of $\mathbb{N}$, there exists a color class that is both an additive and a multiplicative central set. Consequently, it satisfies both versions of the Central Sets Theorem. This naturally leads to the question: 
\textit{how do these two Central Sets Theorems interact?} In this draft, we explore the interplay between the additive and multiplicative forms of the Central Sets Theorem.

\begin{theorem}[\textbf{Unified Central Sets Theorem}]\label{central commutative}
For every finite coloring of $\mathbb{N},$ there exists a color class $A \subseteq \mathbb{N}$ and two functions 
\[
\alpha, \beta: \mathcal{P}_{f}\left(^{\mathbb{N}}\mathbb{N}\right)\rightarrow \mathbb{N} \quad \text{and} \quad H: \mathcal{P}_f(^{\mathbb{N}}\mathbb{N}) \to \mathcal{P}_f(\mathbb{N})
\]
such that:

\begin{itemize}
    \item[(1)] If $F, G \in \mathcal{P}_{f}\left(^{\mathbb{N}}\mathbb{N}\right)$ and $F \subsetneq G$, then 
    \[
    \max H(F) < \min H(G).
    \]

    \item[(2)] For any $m \in \mathbb{N}$ and $G_1 \subsetneq G_2 \subsetneq \dots \subsetneq G_m$ in $\mathcal{P}_{f}\left(^{\mathbb{N}}\mathbb{N}\right)$, and for each $i = 1, 2, \dots, m$ with $f_i \in G_i$, we have:
    
    \begin{itemize}
        \item[(a)] 
        \[
        \sum_{i=1}^{m} \left( \alpha(G_i) + \sum_{t \in H(G_i)} f_i(t) \right) \in A.
        \]
        
        \item[(b)] 
        \[
        \prod_{i=1}^{m} \left( \beta(G_i) \cdot \prod_{t \in H(G_i)} f_i(t) \right) \in A.
        \]
        
        \item[(c)] For all $N \leq m$, 
        \[
        \left( \sum_{i=1}^{N} \left( \alpha(G_i) + \sum_{t \in H(G_i)} f_i(t) \right) \right) 
        \cdot 
        \left( \prod_{i=N+1}^{m} \left( \beta(G_i) \cdot \prod_{t \in H(G_i)} f_i(t) \right) \right) \in A.
        \]
    \end{itemize}
\end{itemize}
\end{theorem}

\section{Preliminaries}
Before proceeding to the proofs of our main results, we recall some essential facts about ultrafilters. For a comprehensive reference, see \cite{HS10}.

Let $\beta \mathbb{N}$ denote the Stone--\v{C}ech compactification of the natural numbers $\mathbb{N}$, i.e., the set of all ultrafilters on $\mathbb{N}$. For any two ultrafilters $p, q \in \beta \mathbb{N}$, define the product $p \cdot q$ by the rule:
\[
A \in p \cdot q \iff \{x \in \mathbb{N} : x^{-1}A \in q\} \in p,
\]
where $x^{-1}A = \{y \in \mathbb{N} : xy \in A\}$. Under this operation, $(\beta \mathbb{N}, \cdot)$ becomes a compact right topological semigroup—that is, the multiplication is associative and continuous in the right variable.

Using Zorn's Lemma, one can show that every compact subsemigroup of $\beta \mathbb{N}$ contains an idempotent element. In particular, both semigroups $(\beta \mathbb{N}, +)$ and $(\beta \mathbb{N}, \cdot)$ possess idempotents.

If $p$ is an idempotent ultrafilter, then any set $A \in p$ is called an \emph{IP-set}. Depending on the semigroup structure, such sets are termed \emph{additive IP-sets} (with respect to $(\beta \mathbb{N}, +)$) or \emph{multiplicative IP-sets} (with respect to $(\beta \mathbb{N}, \cdot)$).

Let $K(\beta \mathbb{N}, +)$ (resp. $K(\beta \mathbb{N}, \cdot)$) denote the smallest two-sided ideal of $(\beta \mathbb{N}, +)$ (resp. $(\beta \mathbb{N}, \cdot)$). The idempotent elements in $K(\beta \mathbb{N}, +)$ (resp. $K(\beta \mathbb{N}, \cdot)$) are known as \emph{additive idempotent ultrafilters} (resp. \emph{multiplicative idempotent ultrafilters}), and the sets they contain are called \emph{additive Central sets} (resp. \emph{multiplicative Central sets}). Similarly, members of the ultrafilters in $K(\beta \mathbb{N}, +)$ (resp. $K(\beta \mathbb{N}, \cdot)$) are known as \emph{additive piecewise syndetic} (resp. \emph{multiplicative piecewise syndetic}) sets.

It is a routine exercise to verify that $\overline{E\left(K\left(\beta \mathbb{N},+\right)\right)}$ is a left ideal of the semigroup $(\beta \mathbb{N},\cdot)$. Hence, the following lemma is immediate.

\begin{lemma}
The set $\overline{E\left(K\left(\beta \mathbb{N},+\right)\right)} \cap E\left(K\left(\beta \mathbb{N},\cdot\right)\right)$ is nonempty.
\end{lemma}

\subsection{Partial Semigroups}

A \emph{partial semigroup} is a pair \((S, \cdot)\), where ``\(\cdot\)'' is a partially defined binary operation on \(S\); that is, it maps a subset of \(S \times S\) into \(S\). It satisfies the following associativity condition: for all \(a, b, c \in S\), if either \((a \cdot b) \cdot c\) or \(a \cdot (b \cdot c)\) is defined, then so is the other, and they are equal.

\begin{ex}
\leavevmode
\begin{itemize}
    \item[(1)]\label{e} Let
    \[
    \mathcal{R} = \left\{ A : \text{there exist } m,n \in \mathbb{N} \text{ such that } A \text{ is an } m \times n \text{ matrix with entries from } \mathbb{Z} \right\},
    \]
    with the usual matrix multiplication. For an \(m \times n\) matrix \(M\) and an \(m' \times n'\) matrix \(N\), define
    \[
    M \cdot N =
    \begin{cases}
        MN & \text{if } n = m', \\
        \text{undefined} & \text{otherwise}.
    \end{cases}
    \]
    Then \((\mathcal{R}, \cdot)\) is a partial semigroup.

    \item[(2)] Let \(\langle x_n \rangle_{n=1}^{\infty}\) be a sequence in a semigroup \((S, \cdot)\), and let
    \[
    FP\left(\langle x_n \rangle\right) = \left\{ \prod_{n \in F} x_n : F \in \mathcal{P}_f(\mathbb{N}) \right\},
    \]
    where the product is taken in increasing order of indices. In general, this set is not closed under the semigroup operation. However, define
    \[
    \left( \prod_{n \in F} x_n \right) \cdot \left( \prod_{n \in G} x_n \right) =
    \begin{cases}
        \prod_{n \in F \cup G} x_n & \text{if } \max F < \min G, \\
        \text{undefined} & \text{otherwise}.
    \end{cases}
    \]
    Then \((T, \cdot)\), where \(T = FP(\langle x_n \rangle)\), is a partial semigroup.
\end{itemize}
\end{ex}

The study of partial semigroups plays an important role in Ramsey theory. For more details, see \cite{HFM}. In this paper, we focus on a special class of partial semigroups known as \emph{adequate partial semigroups}, which are particularly useful for combinatorial applications.

\begin{definition}[\textbf{Adequate Partial Semigroup}]
Let \((S, \cdot)\) be a partial semigroup.
\begin{itemize}
    \item[(1)] For \(s \in S\), define \(\phi(s) = \{ t \in S : s \cdot t \text{ is defined} \}\).
    
    \item[(2)] For \(H \in \mathcal{P}_f(S)\), define \(\sigma(H) = \bigcap_{s \in H} \phi(s)\).
    
    \item[(3)] The partial semigroup \((S, \cdot)\) is said to be \emph{adequate} if \(\sigma(H) \neq \emptyset\) for all finite \(H \subseteq S\).
\end{itemize}
\end{definition}
In Example \ref{e}, one can verify that unlike \((\mathcal{R}, \cdot)\), the partial semigroup \((T, \cdot)\) is \emph{adequate}.

In the case of \((\mathcal{R}, \cdot)\), observe that for any finite subset \(H \subseteq \mathcal{R}\), we have \(\sigma(H) \neq \emptyset\) if and only if all matrices in \(H\) have the same number of columns. That is, if we define
\[
\mathcal{H} = \left\{ A \in \mathcal{R} : A \text{ is a matrix with } r \text{ columns for some fixed } r \in \mathbb{N} \right\},
\]
then \(\sigma(H) \neq \emptyset\) if and only if \(H \subseteq \mathcal{H}\) for some fixed \(r\). Therefore, \((\mathcal{R}, \cdot)\) is not adequate in general.

\subsubsection{Algebra of the Stone-\v{C}ech Compactification of Discrete Partial Semigroups}

Let $(S,\cdot)$ be a partial semigroup, and let $\beta S$ denote the set of ultrafilters on $S$. Then $(\beta S,\cdot)$ is also a partial semigroup.

One of the key advantages of working with \emph{adequate partial semigroups} is that, unlike general partial semigroups, we can identify a genuine semigroup structure within $\beta S$.

\begin{definition}
Given a partial semigroup $(S,\cdot)$, define:
\[
\delta S = \bigcap_{x \in S} \overline{\phi(x)} = \bigcap_{H \in \mathcal{P}_f(S)} \overline{\sigma(H)}.
\]
\end{definition}

Clearly $\delta S \subseteq \beta S$. Remarkably, this set $\delta S$, equipped with a suitably defined operation, forms a semigroup.

For $x \in S$ and $A \subseteq S$, define:
\[
x^{-1}A = \{ y \in \phi(x) : x \cdot y \in A \}.
\]

The following lemma describes the algebraic structure of adequate partial semigroups within the Stone–\v{C}ech compactification.

\begin{lemma}\textup{\cite[Lemma 2.4]{HM}}
Let $S$ be an adequate partial semigroup.
\begin{itemize}
    \item[(a)] Let $x \in S$, $q \in \overline{\phi(x)}$, and $A \subseteq S$. Then
    \[
    A \in x \cdot q \iff x^{-1}A \in q.
    \]
    \item[(b)] Let $p \in \beta S$, $q \in \delta S$, and $A \subseteq S$. Then
    \[
    A \in p \cdot q \iff \{ x \in S : x^{-1}A \in q \} \in p.
    \]
\end{itemize}
\end{lemma}

The following theorem guarantees the existence of idempotents in partial semigroups.
\begin{theorem}\textup{\cite[Theorem 2.6]{HM}}
Let $S$ be an adequate partial semigroup. Then, with the operation described above, $\delta S$ is a compact right topological semigroup and hence contains idempotents.
\end{theorem}

\section{Proof of Theorem \ref{central commutative}}

Let $S$ be a semigroup. For $m \in \mathbb{N}$, define
\[
\mathcal{I}_m = \left\{ (H_1, H_2, \ldots, H_m) : 
\begin{array}{l}
    \text{each } H_j \in \mathcal{P}_f(\mathbb{N}) \text{ and} \\
    \text{for all } j \in \{1, 2, \ldots, m-1\},\ \max H_j < \max H_{j+1}
\end{array}
\right\}.
\]

Given $m \in \mathbb{N}$, $a \in S^{m+1}$, $H \in \mathcal{I}_m$, and $f \in {}^\mathbb{N} S$, define
\[
x(m, a, H, f) = \left( \prod_{j=1}^m \left( a(j) \cdot \prod_{t \in H(j)} f(t) \right) \right) \cdot a(m+1).
\]

Let $\mathcal{I} = \bigcup_{m=1}^\infty \mathcal{I}_m$, and define a partial binary operation $*$ on $\mathcal{I}$ by
\[
(H_1, \ldots, H_m) * (K_1, \ldots, K_n) = (H_1, \ldots, H_m, K_1, \ldots, K_n),
\]
whenever $\max H_m < \min K_1$; otherwise, the operation is undefined.

It follows from \cite{HS10} that $(\mathcal{I}, *)$ is an adequate partial semigroup, and
\[
\delta \mathcal{I} = \bigcap_{n=1}^\infty \overline{ \{ H \in \mathcal{I} : \min H_1 > n \} }^{\delta \mathcal{I}}.
\]

We need the following lemma to prove our main result.

\begin{lemma}[\cite{HS10}, \cite{HS12}, Lemmas 2.10, 14.8.2, 14.8.3]\label{lem:key}
Let $(S, \cdot)$ be a semigroup, $A \subseteq S$ a piecewise syndetic set, and $F \in \mathcal{P}_f({}^\mathbb{N} S)$. Define
\[
\mathcal{B}(A, F, \cdot) = \left\{ (H_1, \ldots, H_m) \in \mathcal{I} : \exists a \in S^{m+1} \text{ such that } \forall f \in F,\ x(m, a, H, f) \in A \right\}.
\]
Then for every $p \in E(\delta \mathcal{I})$, we have $\mathcal{B}(A, F, \cdot) \in p$.
\end{lemma}
Now we are in the position to prove our main result.

\begin{proof}[\textbf{Proof of Theorem \ref{central commutative}:}]

Let  
\[
p \in \overline{E\left(K\left(\beta \mathbb{N},+\right)\right)} \cap E\left(K\left(\beta \mathbb{N},\cdot\right)\right),
\]  
and let \( A \in p \). Then \( A \) is both additively and multiplicatively central.

Define
\[
A^{\star} = \{ x \in A : x^{-1}A \in p \}.
\]
Since \( p \cdot p = p \), it follows that \( A^{\star} \in p \). Moreover, by \cite[Lemma 4.14]{HS12}, if \( x \in A^{\star} \), then
\[
x^{-1} A^{\star} \in p.
\]

Now, by the choice of \( p \), we have \( A^{\star} \in q \) for some additive minimal idempotent \( q \in E(K(\beta \mathbb{N}, +)) \).

Define
\[
A^{\star\star} = \{ x \in A^{\star} : -x + A^{\star} \in q \}.
\]
Then \( A^{\star\star} \in q \). Furthermore, for any \( x \in A^{\star\star} \), it follows that
\[
 -x + A^{\star\star} \in p.
\]

We define \( \alpha(F), \beta(F) \in \mathbb{N} \) and \( H(F) \in \mathcal{P}_{f}(\mathbb{N}) \)
for \( F \in \mathcal{P}_{f}({}^{\mathbb{N}}\mathbb{N}) \) by induction on \( |F| \), satisfying the following inductive hypothesis:

\begin{itemize}
    \item[(1)] For \( F, G \in \mathcal{P}_{f}({}^{\mathbb{N}}\mathbb{N}) \) with \( \emptyset \neq G \subsetneq F \), we have
    \[
    \max H(G) < \min H(F).
    \]

    \item[(2)] Whenever \( n \in \mathbb{N} \), \( G_1, G_2, \ldots, G_n \in \mathcal{P}_{f}({}^{\mathbb{N}}\mathbb{N}) \) with
    \[
    G_1 \subsetneq G_2 \subsetneq \cdots \subsetneq G_n = F,
    \]
    and for each \( i \in \{1,2,\ldots,n\} \), \( f_i \in G_i \), then the following hold:

    \begin{itemize}
        \item[(a)] 
        \[
        \sum_{i=1}^{n} \left( \alpha(G_i) + \sum_{t \in H(G_i)} f_i(t) \right) \in A^{\star\star} \subset A^{\star}.
        \]

        \item[(b)] 
        \[
        \prod_{i=1}^{n} \left( \beta(G_i) \prod_{t \in H(G_i)} f_i(t) \right) \in A^{\star}.
        \]

        \item[(c)] For any \( N \leq n \),
        \[
        \left( \sum_{i=1}^{N} \left( \alpha(G_i) + \sum_{t \in H(G_i)} f_i(t) \right) \right)
        \left( \prod_{i=N+1}^{n} \left( \beta(G_i) \prod_{t \in H(G_i)} f_i(t) \right) \right) \in A^{\star}.
        \]
    \end{itemize}
\end{itemize}

Let \( F = \{ f \} \). Since \( A^{\star\star} \) is an additively piecewise syndetic set, it follows from Lemma~\ref{lem:key} that there exist \( a \in \mathbb{N} \) and \( L_0 \in \mathcal{P}_f(\mathbb{N}) \) such that  
\[
a + \sum_{t \in L_0} f(t) \in A^{\star\star}.
\]

Define  
\[
B = A^{\star\star} \cap \left( -\left( a + \sum_{t \in L_0} f(t) \right) + A^{\star\star} \right) \in q
\quad \text{and} \quad
C = A^{\star} \cap \left( a + \sum_{t \in L_0} f(t) \right)^{-1} A^{\star} \in p.
\]

Using Lemma~\ref{lem:key} again, and noting that \( B \) and \( C \) are additively and multiplicatively piecewise syndetic sets respectively, we have  
\[
\mathcal{B}(B, \{f\}, +) \cap \mathcal{B}(C, \{f\}, \cdot) \neq \emptyset.
\]
In particular, there exist \( a_1, b_1 \in \mathbb{N} \) and \( L_1 \in \mathcal{P}_f(\mathbb{N}) \) such that  
\[
a_1 + \sum_{t \in L_1} f(t) \in B
\quad \text{and} \quad
b_1 \cdot \prod_{t \in L_1} f(t) \in C.
\]
Consequently, we obtain:
\begin{itemize}
    \item \( a_0 + \sum_{t \in L_0} f(t), \; a_1 + \sum_{t \in L_1} f(t) \in A^{\star} \),
    \item \( b_1 \cdot \prod_{t \in L_1} f(t) \in A^{\star} \),
    \item \( \left( a_0 + \sum_{t \in L_0} f(t) \right) \cdot \left( b_1 \cdot \prod_{t \in L_1} f(t) \right) \in A^{\star} \).
\end{itemize}

Define \(\alpha\left(\left\{ f\right\} \right) = a_{1}\) and \(H\left(\left\{ f\right\} \right) = L_{1}\).  
Then both inductive hypotheses are satisfied in this base case.

Now assume \(|F| > 1\), and suppose that for all proper subsets \(G \subsetneq F\),  
the values \(\alpha(G) \in \mathbb{N}\) and \(H(G) \in \mathcal{P}_f(\mathbb{N})\) have already been defined  
satisfying the inductive conditions.

Let  
\[
K = \bigcup \left\{ H(G) : \emptyset \neq G \subsetneq F \right\},
\quad\text{and set}\quad
m = \max K.
\]
Define the sets:
\[
M_{1} = \left\{
\sum_{i=1}^{n} \left( \alpha(G_i) + \sum_{t \in H(G_i)} f_i(t) \right) :
\begin{array}{l}
n \in \mathbb{N},\,
G_1 \subsetneq G_2 \subsetneq \dots \subsetneq G_n \subsetneq F, \\
\langle f_i \rangle_{i=1}^{n} \in \times_{i=1}^{n} G_i
\end{array}
\right\},
\]
\[
M_{2} = \left\{
\prod_{i=1}^{n} \left( \beta(G_i) \cdot \prod_{t \in H(G_i)} f_i(t) \right) :
\begin{array}{l}
n \in \mathbb{N},\,
G_1 \subsetneq G_2 \subsetneq \dots \subsetneq G_n \subsetneq F, \\
\langle f_i \rangle_{i=1}^{n} \in \times_{i=1}^{n} G_i
\end{array}
\right\}.
\]

and 
\begin{align*}
	M_{3} = \Bigg\{  \Bigg( \sum_{i=1}^{N}  \Bigg( \alpha\left(G_{i}\right) + \sum_{t\in H\left(G_{i}\right)} f_{i}\left(t\right) \Bigg)  \Bigg) 
	\cdot \Bigg( \prod_{i=N+1}^{n}  \Bigg( \beta\left(G_{i}\right) \prod_{t\in H\left(G_{i}\right)} f_{i}\left(t\right) \Bigg)  \Bigg) : \\
	\qquad\quad n, N\in\mathbb{N} \text{ with } N\leq n,\ 
	G_{1}\subsetneq G_{2}\subsetneq\ldots\subsetneq G_{n}\subsetneq F,\ 
	\langle f_{i}\rangle_{i=1}^{n}\in\times_{i=1}^{n} G_{i} 
	\Bigg\}.
\end{align*}

Let
$$
B_{1}=A^{\star\star}\cap\bigcap_{x\in M_{1}}\left(-x+A^{\star\star}\right).
$$ 
and
$$
C_{1}=A^{\star}\cap\bigcap_{x\in M_{1}}\left(x^{-1}A^{\star}\right)\bigcap_{x\in M_{2}}\left(x^{-1}A^{\star}\right)\bigcap_{x\in M_{3}}\left(x^{-1}A^{\star}\right).
$$
 By hypothesis (2), $M_1 \subset A^{\star\star} \subset A^\star$ and $M_2 \subset A^\star$, which, together with the fact that the sets $M_i$ for $i = 1, 2, 3$ are finite, implies that 
\[
B_1 \in q \quad \text{and} \quad C_1 \in p.
\]
Thus, Lemma \ref{lem:key} similarly guarantees that there exist
$a, b \in \mathbb{N}$ and $L \in \mathcal{P}_{f}(\mathbb{N})$ such that $\min L > m$, and for all $f \in F$,
\[
a + \sum_{t \in L} f(t) \in B_1 \quad \text{and} \quad b \cdot \prod_{t \in L} f(t) \in C_1.
\]

Let $\alpha(F) = a$, $\beta(F) = b$, and $H(F) = L$. 
Since $\min L > m$, the first inductive hypothesis is satisfied. 
We now verify hypothesis (2). If $n = 1$, then we have
\[
a + \sum_{t \in L} f(t) \in B_1 \subset A^{\star\star} \subseteq A^{\star}
\quad \text{and} \quad 
b \cdot \prod_{t \in L} f(t) \in C_1 \subseteq A^{\star},
\]
which shows that hypothesis (2) holds in this case.

And if $n > 1$, let
\[
x = \sum_{i=1}^{n-1} \left( \alpha(G_i) + \sum_{t \in H(G_i)} f_i(t) \right), \quad
y = \prod_{i=1}^{n-1} \left( \beta(G_i) \prod_{t \in H(G_i)} f_i(t) \right),
\]
and
\[
z = \left( \sum_{i=1}^{N} \left( \alpha(G_i) + \sum_{t \in H(G_i)} f_i(t) \right) \right)
    \left( \prod_{i=N+1}^{n} \left( \beta(G_i) \prod_{t \in H(G_i)} f_i(t) \right) \right).
\]

 Then $x \in M_1$, $y \in M_2$ and $z \in M_3$, and hence
\begin{itemize}
    \item $a + \sum_{t \in L} f(t) \in B_1 \subseteq A^{\star\star} - y \subseteq A^{\star} - y$;
    
    \item $b \cdot \prod_{t \in L} f(t) \in C_1 \subseteq z^{-1} A^{\star}$;
    
    \item For $N \leq n$, $b \cdot \prod_{t \in L} f(t) \in C_1 \subset z^{-1} A^{\star}$
\end{itemize}
for all $G_1 \subsetneq G_2 \subsetneq \ldots \subsetneq G_n = F$.

Therefore, hypothesis (2) holds for all $n \in \mathbb{N}$.  
This completes the inductive proof of the central set theorem.

\end{proof}

 \section*{Acknowledgement}  The second author of this paper is supported by NBHM postdoctoral fellowship with reference no: 0204/27/(27)/2023/R \& D-II/11927.

\bibliographystyle{plain}

\end{document}